\magnification=\magstep1 
\overfullrule=0pt
\def\eqde{\,{\buildrel \rm def \over =}\,} 
    
       \def\i{{\rm i}} 
 \def\eps{\epsilon}    
         
  \def\L{{\cal L}}

\font\huge=cmr10 scaled \magstep2

\def\z{{\cal Z}}  
   
\input amssym.def
\def\Z{{\Bbb Z}}  \def\Q{{\Bbb Q}}  
\def\C{{\Bbb C}}

\rightline{ESI preprint \# 611}
\rightline{math.QA/9902064}
\rightline{January, 1999}
\bigskip\medskip
\centerline{{\bf \huge  The Cappelli-Itzykson-Zuber  A-D-E Classification}}
\bigskip \bigskip   \centerline{Terry Gannon}\medskip
\centerline{{\it Department of Mathematical Sciences, University of Alberta,}}
\centerline{{\it Edmonton, Alberta, Canada, T6G 2G1} } \smallskip
\centerline{{e-mail: tgannon@math.ualberta.ca}}
\bigskip\medskip

In 1986 Cappelli, Itzykson and Zuber classified all modular invariant partition
functions for the conformal field theories associated to the affine $A_1$
algebra; they
found they fall into an A-D-E pattern. Their proof was difficult and
attempts to generalise it to the other affine algebras failed -- in
hindsight the reason is that their argument ignored most of the rich 
mathematical structure
present. We give here the ``modern'' proof of their result; it is an
order of magnitude simpler and shorter, and much of it has already been
extended to all other affine algebras. We conclude with some remarks on the 
A-D-E pattern appearing in this and other RCFT classifications.
\bigskip\bigskip\bigskip\bigskip

\noindent{{\bf 1. The problem}}\bigskip 

 One of the more important results in conformal field theory is
surely the classification due to Cappelli, Itzykson, and Zuber [1; see also
2] of the
genus 1 partition functions for the theories associated to $A_1^{(1)}$
(which in turn implies the classification of the minimal models).
Their list was curious: the partition functions fall 
into the A-D-E
pattern familiar from the simply-laced Lie algebras,
finite subgroups of SU$_2(\C)$, simple singularities,
 subfactors with index $<4$, representations of quivers, etc. See e.g.\ [3].

The problem can be phrased as follows. Fix any integer $n\ge 3$. Let
$P_+=\{1,2,\ldots,n-1\}$, and let $S$ and $T$ be the $(n-1)\times (n-1)$
matrices with entries
$$S_{ab}=\sqrt{{2\over n}}\,\sin(\pi\,{ab\over n})\ ,\qquad
\qquad T_{ab}=\exp[\pi\i\,{a^2\over 2n}]\,\delta_{a,b}\ .$$
Find all $(n-1)\times (n-1)$ matrices $M$ such that
\smallskip\item{$\bullet$} $M$ commutes with $S$ and $T$: $MS=SM$ and $MT=TM$

\item{$\bullet$} $M$ has nonnegative integer entries: $M_{ab}\in \Z_+$ for all $a,b\in P_+$

\item{$\bullet$} $M$ is normalised so that $M_{11}=1$\ .\smallskip
\noindent Call any such $M$ a {\it physical invariant}. Since most entries
$M_{ab}$ are usually zero, it is more convenient
to formally express $M$ as the coefficient matrix for the combination
$$\z=\sum_{a,b=1}^{n-1} M_{ab}\,\chi_a\,\chi_b^*\ .$$

{{\bf Theorem}} [1]. {\it The complete list of
physical invariants is (using $Ja\eqde n-a$)}
$$\eqalignno{{\cal A}_{n-1}=&\,\sum_{a=1}^{n-1}\,|\chi_a|^2\ ,\qquad\qquad
\qquad \forall n\ge 3&
\cr{\cal D}_{{n\over 2}+1}=&\,\sum_{a=1}^{n-1}\,\chi_a\,\chi_{J^aa}^*\ ,\qquad 
\qquad\qquad{\rm whenever}\ {n\over 2}\ {\rm is\ even}&\cr
{\cal D}_{{n\over 2}+1}=&\,|\chi_1+\chi_{J1}|^2+|\chi_3+\chi_{J3}|^2+\cdots
+2|\chi_{{n\over 2}}|^2\ ,\qquad {\rm whenever}\ {n\over 2}\ {\rm is\ odd}&\cr
{\cal E}_6=&\,|\chi_1+\chi_7|^2+|\chi_4+\chi_8|^2+|\chi_5+\chi_{11}|^2\ ,
\qquad\qquad {\rm for}\ n=12&\cr
{\cal E}_7=&\,|\chi_1+\chi_{17}|^2+|\chi_5+\chi_{13}|^2+|\chi_7+\chi_{11}|^2
&\cr&\,+\chi_9\,(\chi_3+\chi_{15})^*+(\chi_3+\chi_{15})\,\chi^*_9+|\chi_9|^2
\ ,\qquad\qquad {\rm for}\ n=18&\cr
{\cal E}_8=&\,|\chi_1+\chi_{11}+\chi_{19}+\chi_{29}|^2+|\chi_7+\chi_{13}+
\chi_{17}+\chi_{23}|^2\ ,\qquad {\rm for}\ n=30\ .&\cr}$$\smallskip

These realise the A-D-E pattern, in the following sense.
The Coxeter number $h$ of the name ${\cal X}_\ell$ equals the
corresponding value of $n$, and the exponents of ${ X}_{\ell}$
(i.e.\ the $m_i$ in the eigenvalues $4\sin^2(\pi\,{m_i\over 2n})$ of its 
Cartan matrix) equal those
$a\in P_+$ for which $M_{aa}\ne 0$. An interpretation of the nondiagonal
entries of $M$ has been provided by Ocneanu [4] (see also [5]).

Cappelli-Itzykson-Zuber proved their theorem by first finding an explicit basis for the 
space of all matrices commuting with $S$ and $T$.
Unfortunately their proof  was long and formidable. 
Considering all of the structure implicit in the problem, we should
anticipate a much more elementary argument. This is
not merely of academic interest, because there is a natural generalisation of
this problem to all other affine algebras. Several people had tried to extend
the argument of [1] to these larger algebras, but with [6] it became clear
that some other approach was necessary, or the generalisation would never
be achieved. And of course another reason is that the more transparent the 
argument, the better the chance of ultimately understanding the connection with A-D-E.

In this paper we provide a considerably shorter proof of the theorem,
bearing no resemblance to the older arguments. 
Our proof  is an example of the ``modern'' approach to physical invariant 
classifications. See [7] for a summary of the current status of these
classifications for the other affine algebras.

The argument which follows is completely elementary: no knowledge of e.g.\
CFT or Kac-Moody algebras is assumed. It is based on various talks I've given,
 most recently at the Schr\"odinger Institute in Vienna
where I wrote up this paper and who I thank for generous hospitality.
I also thank D.\ Evans, M.\ Flohr, J.\ McKay, V.\ Petkova, and J.-B.\ Zuber
for correspondence.

\vfill\eject
\noindent{{\bf 2. The combinatorial background}} 

\bigskip In this section we include some of the basic tools belonging
to any classification of the sort, and we give a flavour of their proofs.
We will state them for the specific problem given above, but
everything generalises without effort [8].

First note that commutation of $M$ with $T$ implies the selection rule
$$M_{ab}\ne 0\quad\Longrightarrow\quad a^2\equiv b^2\
\quad({\rm mod}\ 4n)\ .\eqno(2.1)$$

Next, let us write down some of the basic properties obeyed by $S$.
$S$ is symmetric and orthogonal (so $M=SMS$), and
$$S_{1b}\ge S_{11}>0\ .\eqno(2.2)$$
The permutation $J$ of $P_+$, defined by $Ja=n-a$, corresponds to the order 2 
symmetry of the extended Dynkin diagram of $A_1^{(1)}$; it satisfies
$$S_{Ja,b}=(-1)^{b+1}S_{ab}\ .\eqno(2.3)$$
Note that the element $1\in P_+$ is both physically and mathematically
special; our strategy will be to find all possible first rows and columns
of $M$, and then for each of these possibilities to find the remaining
entries of $M$.

The easiest result follows by evaluating $MS=SM$ at $(1,a)$ for any $a\in P_+$:
$$S_{1a}+\sum_{b=2}^{n-1}M_{1b}\,S_{ba}\ge 0\ ,\eqno(2.4)$$
with equality iff the $a$th column of $M$ is identically 0. Equation (2.4)
has two uses: it severely constrains the values of $M_{1b}$ (similarly
$M_{b1}$), and it says precisely which columns (and rows) are nonzero.

Another simple observation is
$$1=M_{11}=\sum_{a,b=1}^{n-1} S_{1a}\,M_{ab}\,S_{1b}\ge S_{11}^2\sum_{a,b=1}^{
n-1}M_{ab}\ .$$
This tells us that each entry $M_{ab}$ is bounded above
by ${1\over S_{11}^2}$ (we will use this below). In particular, there can
only be finitely many physical invariants for each $n$. (This same calculation 
shows more generally that there will only be finitely many physical invariants
 for a given affine algebra $X_r^{(1)}$ and level $k$.)

Next, let's apply the triangle inequality to sums involving (2.3).
Choose any $i,j\in\{0,1\}$. Then
$$M_{J^i1,J^j1}=\sum_{a,b=1}^{n-1}(-1)^{(a+1)i}\,S_{1a}\,M_{ab}\,(-1)^{(b+1)j}\,
S_{1b}\ .$$
Taking absolute values, we obtain
$$M_{J^i1,J^j1}\le\sum_{a,b=1}^{n-1}S_{1a}\,M_{ab}\,S_{1b}=M_{11}=1\ .$$
Thus $M_{J^i1,J^j1}$ can equal only 0 or 1. If it equals 1, then we obtain
the selection rule: $(a+1)i\equiv (b+1)j$ (mod 2) whenever $M_{ab}\ne 0$;
this implies the symmetry $M_{J^ia,J^jb}=M_{ab}$ for all $a,b\in P_+$.

Whenever you have nonnegative matrices in your problem, and it makes sense
to multiply those matrices, then you should seriously consider using
Perron-Frobenius theory -- a collection of results concerning
the eigenvalues and eigenvectors of nonnegative matrices. Our
$M$ is nonnegative, and although multiplying $M$'s may not give us back a
physical invariant, at least it will give us a matrix commuting with
$S$ and $T$. In other words, the commutant is much more than merely a
vector space, it is in fact an algebra.

Important applications of this thought are the following two lemmas.

\medskip{\bf Lemma 1}.\quad {\it Let $M$ be a physical invariant,
and suppose $M_{a1}=\delta_{a,1}$ -- i.e.\ the first column of $M$ is all
zeros except for $M_{11}=1$. Then $M$ is a permutation matrix -- i.e.\
there is some permutation $\pi$ of $P_+$ such that $M_{ab}=\delta_{b,\pi a}$,
and $S_{\pi a,\pi b}=S_{ab}$.}
\medskip

This is proved by studying the powers $(M^T
M)^L$ as $L$ goes to infinity: its diagonal entries will grow exponentially
with $L$, unless there is at most one nonzero entry on each column of $M$,
and it equals 1. (Recall that the entries of $(M^TM)^L$ must be bounded
above.) Then show that also $M_{1a}=\delta_{1,a}$
 (evaluate $MS=SM$ at (1,1)), and look at $(M\,M^T)^L$.
 Lemma 1 was found independently by Schellekens and Gannon.

That argument is elementary enough that it required no knowledge of 
Perron-Frobenius. But Perron-Frobenius is needed for generalisations.
In this fancier language, what the preceding argument  shows is: write
$M$ as the direct sum of indecomposable submatrices; then the largest
eigenvalue of the submatrix containing (1,1)
bounds above that for each other submatrix. Arguing with a little
more sophistication, we obtain much more. The special case we need is: 

\medskip{\bf Lemma 2} [8].\quad {\it Let $M$ be a physical 
invariant,
and suppose $M_{a1}\ne 0$ only for $a=1$ and $a=J1$, and similarly for
$M_{1a}$ -- i.e.\ the first row and column of $M$ are all
zeros except for $M_{J^i1,J^j1}=1$. Then the $a$th row (or column) of $M$
will be identically 0 iff $a$ is even. Moreover, let $a,b\in P_+$, both
different from ${n\over 2}$, and suppose $M_{ab}\ne 0$. Then
$$M_{ac}=\left\{\matrix{1&{\rm if}\ c=b\ {\rm or}\ c=Jb\cr 0&{\rm otherwise}
\cr}\right.$$
and a similar formula holds for $M_{cb}$.}
\medskip

This lemma says that the indecomposable submatrices of $M$ which don't
involve ${n\over 2}$ (the fixed-point of $J$) will either be trivial (0)
(for even places on the diagonal), or involve blocks $\left(\matrix{1&1\cr 
1&1}\right)$. You can check this for
the ${\cal D}_{even}$ and ${\cal E}_7$ partition functions. 

Our final ingredient is a Galois symmetry obeyed by $S$, and its consequence
for $M$. Again, see e.g.\ [8] for a proof. Let $\L$ be the set of all $\ell$
coprime to $2n$. For each $\ell\in \L$, there is a permutation $a\mapsto 
[\ell a]$ of $P_+$, and a choice
of signs $\eps_\ell:P_+\rightarrow\{\pm 1\}$, such that
$$M_{ab}=\epsilon_\ell(a)\,\epsilon_\ell(b)\,M_{[\ell a],[\ell b]}\ ,
\eqno(2.5)$$
for all $a,b\in P_+$.
In particular, write $\{x\}$ for the unique number congruent to $x$ (mod $2n$)
satisfying $0\le \{x\}<2n$. Then if $\{\ell a\}<n$, put $[\ell a]=\{\ell a\}$
and $\eps_\ell(a)=+1$, while if $\{\ell a\}>n$, put $[\ell a]=2n-\{\ell a\}$
and $\eps_\ell(a)=-1$. This `Galois symmetry' (2.5) comes from hitting
$M=SMS$ with the $\ell$th `Galois automorphism'. Any polynomial over $\Q$
with a $2n$th root of unity $\zeta$ as a zero -- and the entries of
$M=SMS$ can be interpreted
in that way -- also has $\zeta^\ell$ as a zero.
 We then use $\sin(\pi\,\ell{ab\over n})=\eps_\ell
(a)\,\sin(\pi\,{[\ell a]b\over n})$. From (2.5) and
 the positivity of $M$, we get for all $\ell\in \L$ the 
Galois selection rule
$$M_{ab}\ne 0\quad\Longrightarrow\quad\epsilon_\ell(a)=
\epsilon_\ell(b)\ .\eqno(2.6)$$
(2.5) and (2.6), valid for any affine algebras, were first found independently
by Gannon and Ruelle-Thiran-Weyers. The Galois interpretation, and extension
to all RCFT, is due to Coste-Gannon.

\bigskip\bigskip\noindent{\bf 3. The ``modern'' proof of the $A_1^{(1)}$
classification}\bigskip

The last section reviewed the basic tools shared by
 all modular invariant partition function
classifications. In this section we specialise to $A_1^{(1)}$.

The first step will be to find all possible values of $a$ such that
$M_{1a}\ne 0$ or $M_{a1}\ne 0$. These $a$ are severely constrained.
We know two generic possibilities: $a=1$ (good for all $n$), and $a=J1$
(good when ${n\over 2}$ is odd). We now ask the question, what other 
possibilities for $a$ are there? Our goal is to prove (3.4). Assume $a\ne 1,
J1$.

There are only two constraints on $a$ which we will need. One is (2.1):
$$(a-1)\,(a+1)\equiv 0\qquad ({\rm mod}\ 4n)\ .\eqno(3.1)$$
More useful is the Galois selection rule (2.6), which we 
can write as $\sin(\pi\ell {a\over n})\sin(\pi\ell{1\over n})>0$, for all
$\ell\in \L$. But a product
of sines can be rewritten as a difference of cosines, so we get
$$\cos(\pi\,\ell\,{a-1\over n})>\cos(\pi\,\ell\,{a+1\over n})\ .\eqno(3.2)$$
Since $\ell$ obeys (3.2) iff $\ell+n$ does, we can take $\ell$ in (3.2) to be
coprime merely to $n$ instead of $2n$. Call $\L'$ the set of these $\ell$.
(3.2)  is strong and easy to solve; here is my argument.

Define $d={\rm gcd}(a- 1,\,2n)$, $d'={\rm gcd}(a+1,\,2n)$. Note from (3.1)
that gcd($
d,d')=2$ and $dd'=4n$, so $d,d'\ge 6$. We can choose $\ell_0,\ell'\in \L'$ so 
that $\ell'\,(a+1)\equiv d'$ (mod $2n$) and
$$\ell_0\,(a-1)\equiv \left\{\matrix{n-d&{\rm if}\ {d\over 2}\ {\rm
is\ odd\ and\ }{n\over 2}\ {\rm is\ even}\cr n-2d&{\rm if}\ {d\over 2}\
{\rm is\ odd\ and\ } {n\over 2}\ {\rm is\ odd}&\cr n-{d\over 2}&{\rm
otherwise,\ i.e.\ if}\ {d'\over 2}\ {\rm is\ odd}\cr}\right. \qquad\qquad
{\rm (mod}\ 2n)\ .$$
Now define $\ell_i={2ni\over d}+\ell_0$. Then $\ell_i
\,(a-1)\equiv \ell_0\,(a-1)$ (mod $2n$) for all $i$, and for $0\le i<{d\over 
2}$ the numbers $\ell_i\,(a+1)$ will all be distinct (mod $2n$).
For those $i$, precisely $\phi({d\over 2})$ of the $\ell_i$
will be in $\L'$, where $\phi(x)$ is the Euler totient, i.e.\ the number
of positive integers less than $x$ coprime to $x$. (This count follows from
the fact that for any prime $p$ dividing ${d\over 2}$, $p$ won't divide
${2n\over d}$ and hence exactly one value of $i$ (mod $p$) will be
forbidden.)

Now, the numbers $\ell_i\,(a+ 1)$ are all multiples of $d'$.
So (3.2) with $\ell=\ell_i$ gives us 
$$(\phi({d\over 2})-1)\,d'<\left\{\matrix{2d&{\rm if}\ {d\over 2}\ {\rm
is\ odd\ and\ }{n\over 2}\ {\rm is\ even}\cr 4d&{\rm if}\ {d\over 2}\
{\rm is\ odd\ and\ } {n\over 2}\ {\rm is\ odd}&\cr d&{\rm otherwise}\cr}
\right. \ .\eqno(3.3)$$
Also, (3.2) with $\ell=\ell'$ requires $d'>d(\ge 6)$.
Combining this with (3.3), we get $\phi({d\over 2})-1<
2$, 4, or 1, which has the solutions $d=6$ (for $n$ some multiple of 4),
and $d=6$ or 10 (for $n$ an odd multiple of 2). (3.3)
now gives us exactly 3 possibilities: $d=6$, $d'=8$, $n=12$ (which 
yields ${\cal E}_6$ as we will see below); $d=6$, $d'=20$,
$n=30$, and $d=10$, $d'=12$, $n=30$ (both which correspond to ${\cal
E}_8$).

So what we have shown is that, provided $n\ne 12,30$, $M$ obeys the
strong condition 
$$M_{a1}\ne 0\ {\rm or}\ M_{1a}\ne 0\qquad\Longrightarrow\qquad a\in\{1,J1\}
\ .\eqno(3.4)$$
Consider first {\bf case 1}: $M_{a1}=\delta_{a,1}$. This is the condition
in Lemma 1, and so we know $M_{ab}= \delta_{
b,\pi a}$ for some permutation $\pi$ of $P_+$ obeying $S_{ab}=S_{\pi a,\pi b}$.
We know $\pi 1=1$; put $m:=\pi 2$. Then $\sin(\pi {2\over n})=\sin(\pi\,
{m\over n})$, and so we get either $m=2$ or $m=J2$. By $T$-invariance (2.1), 
the second possibility can only occur if $4\equiv (n-2)^2$ (mod $4n$), i.e.\
4 divides $n$. But for those $n$ ${\cal D}_{{n\over 2}+1}$ is also a permutation
matrix, so replacing $M$ if necessary with the matrix product $M\,{\cal 
D}_{{n\over 2}+1}$, we can always require $m=2$, i.e.\ $\pi 2=2$.

Now take any $a\in P_+$ and write $b=\pi a$: we have both $\sin(\pi\,{a\over
n})=\sin(\pi\,{b\over n})$ and $\sin(\pi\,{2a\over n})=\sin(\pi\,{2b\over n})$.
Dividing these gives $\cos(\pi\,{a\over n})=\cos(\pi\,{b\over n})$, and
we read off that $b=a$, i.e.\ that $M$ is the identity matrix ${\cal A}_{n-1}$.

The other possibility, {\bf case 2}, is that both $M_{1,J1}\ne 0$ and $M_{J1,
1}\ne 0$. Then Lemma 2 applies. (2.1) says $1\equiv (n-1)^2$ (mod $4n$),
i.e.\ ${n\over 2}$ is odd. $n=6$ is trivial (the only unknown entry,
$M_{3,3}$, is fixed by $MS=SM$ at (1,3)),
so consider $n\ge 10$. The role of `2' in {\bf case 1} will be played here
by `3'. The only difference is the complication caused by the fixed-point
${n\over 2}$. Can $M_{3,{n\over 2}}\ne 0$? If so, then Lemma 2 would imply
$M_{3,a}=0$ for all $a\ne{n\over 2}$. Evaluating $MS=SM$ at $(3,1)$,
we obtain $M_{3,{n\over 2}}= 2\sin(\pi\,{3\over n})$, i.e.\ $n= 18$,
which corresponds to ${\cal E}_7$ as we show later.

Thus we can assume for now that both $M_{3,{n\over 2}}=M_{{n\over 2},3}=0$, 
and so by
Lemma 2 there will be a unique $m<{n\over 2}$ for which $M_{3,m}\ne 0$.
$MS=SM$ at $(3,1)$ now gives $m=3$. For any odd $a\in P_+$, $a\ne {n\over
2}$, can we have $M_{{n\over 2},a}\ne 0$? 
If so then  $MS=SM$ at $(1,a)$ and $(3,a)$ 
would give us $2\sin(\pi\,{a\over n})=M_{{n\over 2},a}=
2\sin(\pi\,{3a\over n})$, which is impossible.
Therefore Lemma 2 again applies, and we get a unique $b<{n\over 2}$ for which
$M_{ba}\ne 0$. The usual argument forces $b=a$, and we obtain the desired
result: $M={\cal D}_{{n\over 2}+1}$.

\vfill\eject\noindent{{\it 3.1. The exceptional at $n=12$}}\smallskip

We know $M_{1a}\ge 1$ for some $a\in P_+$ with gcd$(a+1,24)=8$ -- i.e.\ $a=7$.
From (2.4) at $a=2$, we get $\sin({\pi\over 6})-M_{1,7}\sin({\pi
\over 6})\ge 0$. Thus $M_{1,7}=1$. Applying the Galois symmetry (2.5) for 
$\ell=
5,7,11$, we obtain the terms $|\chi_1+\chi_7|^2+|\chi_5+\chi_{11}|^2$
in ${\cal E}_6$. Now use (2.4) to show that among the remaining entries of $M$,
only the 4th and 8th rows and columns will be nonzero. $M_{J1,J1}=1$ tells
us $M_{44}=M_{88}$ and $M_{84}=M_{48}$. These must be equal, by evaluating
$MS=SM$ at (4,2), and now either Perron-Frobenius or $MS=SM$ at
$(1,4)$ forces that common value to be 1. We thus obtain $M={\cal E}_6$.

\medskip\noindent{{\it 3.2. The exceptional at $n=18$}}\smallskip

We know $M_{3,9}=1$ and that $M_{3,a}=0$ for all other $a\ne 9$.
$T$-invariance (2.1) and Lemma 2 applied to the other odd $a<9$,  force
$M_{aa}=1$. The only remaining entry is $M_{9,9}$,
which is fixed by $MS=SM$ at (9,1). We get $M={\cal E}_7$.

\medskip\noindent{{\it 3.3. The exceptional at $n=30$}}\smallskip

We know either $M_{1,11}$ or
$M_{1,19}$ is nonzero; the only other (potentially) nonzero $M_{1a}$ are at
$a=1,J1$. Suppose first that $M_{1,J1}=1$, so
$M_{1,11}=M_{1,19}$. Then (2.4) at $a=3$ forces $M_{1,11}=1$; Galois (2.5) for
$\ell=7,11,13,17,19,23,29$ gives us rows 7, 11, 13, 17, 19, 23, 29 of $M$,
 and (2.4) tells us all other rows must vanish.

If instead $M_{1,J1}=0$, then (2.4) at $a=2,3,4$ gives our contradiction.

\bigskip\bigskip\noindent{\bf 4. Closing remarks}\bigskip

There are two reasons to be optimistic about the possibilities of a
classification of all modular invariant partition functions ($=\,$physical 
invariants) for all simple $X_r$. One is the main general result in the 
problem [8],
which gives the analogue for any $X_r$ of the $A_1^{(1)}$ physical
invariants named ${\cal A}_{\star}$, ${\cal D}_{\star}$, and
${\cal E}_7$. See [7] for a discussion.
The other cause for optimism is the shortness and simplicity of the
above proof for $A_1^{(1)}$.

The reader should be warned though that $A_1^{(1)}$ is an exceptionally
gentle case -- as we've seen, the proof quickly reduces essentially to combinatorics. 
Our argument here is a
projection of the general argument onto this special case, and this projection
loses most of the structure present in the general proofs. The general
arguments are necessarily more subtle and sophisticated. Nevertheless
this paper should help the interested reader understand the further
literature on this fascinating problem, and  make more accessible the proof of
the important classification of Cappelli-Itzykson-Zuber.

A big question is, does this new proof shed any light on the main mystery here:
the A-D-E pattern to our Theorem? It does not appear
to. But it should be remarked that it is entirely without foundation to argue
 that this $A_1^{(1)}$ classification is `equivalent' to any other A-D-E
 one. There is a connection with the other A-D-E classifications which
 should be explained, and which has not yet been satisfactorily explained. But 
what we
 should look for is some critical combinatorial part of a proof which
 can be identified with critical parts in other A-D-E 
classifications. For instance, does an argument equivalent to that
surrounding (3.2) appear elsewhere in the A-D-E literature?

There has been some progress elsewhere at understanding our
A-D-E. Nahm [9] constructed the invariant ${\cal X}_\ell$ in terms of the
compact simply-connected Lie group of type $X_\ell$, and in this way could
interpret the $n=h$ and $M_{m_im_i}\ne 0$ coincidences.
A very general explanation for A-D-E has been suggested by
Ocneanu [4] using his theory of paragroups and path algebras on
graphs,  but unfortunately  it has
yet to be published (though details of this work are slowly appearing
elsewhere -- see e.g.\ [5]). Related to this is the
work by Di Francesco, Petkova and Zuber on fusion graphs (see e.g.\ [10]);
for its interpretation involving $II_1$ subfactors see e.g.\ [5].
Also worth mentioning is the classification of boundary conditions
in CFT (i.e.\ of partition functions on a finite cylinder rather than
a torus). This appears to be equivalent to a classification of
certain $\Z_+$-valued representations of the fusion ring [11]; for
$A_1^{(1)}$ the problem is quickly reduced to considering 
symmetric $\Z_+$-matrices with largest eigenvalue $<2$, and from this we
once again get an A-D-E pattern.
Nevertheless, the A-D-E in CFT seems to remain almost as mysterious
now as  it did a dozen years ago...

Clearly, a very interesting question is, what form if any does the A-D-E
pattern take for $A_2^{(1)}$ physical invariants? $A_3^{(1)}$? etc.
A step in this direction is provided by the physical invariant
$\leftrightarrow$ fusion graph $\leftrightarrow$ subfactor theory
[5,10] alluded to above. In particular this interprets and generalises
the $M_{m_im_i}\ne 0$ coincidence (at least for the so-called
`block-diagonal' physical invariants, i.e.\ ${\cal Z}$ which can be
expressed as sums of squares: ${\cal Z}=\sum|\chi+\chi'+\cdots|^2$).
Related to this is the following.
It is known that orbifolding a 4-dimensional $N=4$ supersymmetric gauge
theory by any finite subgroup $G\subset {\rm SU}_2({\Bbb C})$ leads to a CFT
with $N=2$ supersymmetry, whose matter matrix (giving numbers of fermions
and scalars) can be read off from the Dynkin diagram corresponding to $G$.
For finite subgroups of SU$_3({\Bbb C})$ and SU$_4({\Bbb C})$, we would
get $N=1$ and $N=0$ supersymmetry, respectively. The (directed) graphs corresponding
to the matter matrices for $G\subset {\rm SU}_3({\Bbb C})$ are given in
[12] and closely resemble the fusion graphs of [10] for $A_2^{(1)}$
physical invariants. Indeed, [12] make the tantalising conjecture that
there exists a McKay-like correspondence between certain singularities of
type ${\Bbb C}^n/G$ (or corresponding orbifold theories) for $G\subset
{\rm SU}_n({\Bbb C})$, and the physical invariants of $A_{n-1}^{(1)}$.
Now, the finite subgroups of SL$_n({\Bbb C})$, at least for $n\le 7$, are
known (Blichfeldt 1917, Brauer 1967, Lindsey 
1971, Wales 1968, for $n=4,5,6,7$ resp.), so presumably the work of [12]
can with effort be extended and their conjecture more precisely stated
and tested. It should be noted though that [12] also makes
use of only those `block-diagonal' $A_2^{(1)}$ physical invariants 
(in analogy with the $A_1^{(1)}$ classification,
it is as if they would ignore ${\cal D}_{odd}$ and ${\cal E}_7$ --
these graphs are also missing from the list of principle graphs of
subfactors).
What if anything should correspond to the remaining physical invariants is unknown.

Incidently, there is a nice little curiousity contained within many
modular invariants: another A-D-E! This
A-D-E applies to {\it any} physical invariant (i.e.\
for any RCFT, not necessarily related to $A_1^{(1)}$) which  looks like
${\cal Z}=|\chi_1+\chi_{1'}|^2+{\rm stuff}$. The label $1'$ can be anything
in $P_+$, and `stuff' can be any sesquilinear combination of $\chi_i$'s,
provided it doesn't contain $\chi_1$ (the vacuum) or $\chi_{1'}$. In other
words, the indecomposable submatrix of $M$ containing $(1,1)$ is required
to be $\left(\matrix{1&1\cr 1&1}\right)$, but otherwise $M$ is unconstrained.
Then to $M$ we can associate several extended Dynkin diagrams of
A-D-E type, as follows.

Put a node on the left of the page for each $a\in P_+$
whose row $M_{a\star}$ is nonzero, and put a node on the right of the
page for each $b\in P_+$ whose column $M_{\star b}$ is nonzero.
Connect $a$ (on the left) and $b$ (on the right) with precisely
$M_{ab}$ edges. The result will be a set of extended Dynkin diagrams
of A-D-E type! (For these purposes we will identify two nodes connected with
2 lines as the extended $A_1$ diagram.)

For example, let's apply this to our $A_1^{(1)}$ classification.
 Any partition function  ${\cal D}_{2\ell}$ is of this kind, and
its corresponding graph will consist of $\ell-1$ diagrams of (extended) 
$A_3$ type, and one of extended $A_1$ type. The exceptional ${\cal E}_6$ consists
of three $A_3$'s, and the exceptional ${\cal E}_7$ consists of
three $A_3$'s and one $D_5$. Again, this fact (proved in [8]) is not 
restricted to the $A_1^{(1)}$ physical invariants.

This little curiousity is not 
as deep or mysterious as the Cappelli-Itzykson-Zuber A-D-E pattern, and has to
do with the $\Z_+$-matrices with largest eigenvalue 2.

There are 4 other claims for A-D-E classifications of families of 
RCFT physical invariants, and all of them inherit their (approximate)
A-D-E pattern from the more fundamental $A_1^{(1)}$ one. The two rigourously
established ones are the $c<1$ minimal models, also proven in [1], and the
$N=1$ superconformal minimal models, proved in [13]. In both cases the
physical invariants are parametrised by pairs of A-D-E diagrams.
The list of known
$c=1$ RCFTs [14] also looks like A-D-E (two series parametrised by $\Q_+$,
and three exceptionals), but the completeness of that list
has never been successfully proved (or at least such a proof has never 
been published).

The fourth classification often quoted as A-D-E, is the $N=2$
superconformal minimal models. The only rigourous classification of
these is accomplished in [15], assuming the generally believed but still
unproven coset realisation $(SU(2)_k\times U(1)_4)/U(1)_{2k+4}$. The 
connection here with A-D-E turns out to
be rather weak: e.g.\ 20, 30, and 24 distinct invariants would 
have an equal right to be  called
${\cal E}_6$, ${\cal E}_7$, and ${\cal E}_8$ respectively. It appears to
this author that the frequent claims that the $N=2$ minimal models fall
into an A-D-E pattern are without serious foundation, or at least require
major reinterpretation.

\vfill\eject \noindent{{\bf References}} \bigskip

\item{1.} A.\ Cappelli, C.\ Itzykson, and J.-B.\ Zuber,  
{\it Commun.\ Math.\ Phys.}\ {\bf 113}, 1 (1987).

\item{2.} A.\ Cappelli, C.\ Itzykson, and J.-B.\ Zuber, {\it Nucl.\ Phys.}\ 
{\bf B280 [FS18]}, 445 (1987);

\item{} D.\ Gepner and Z.\ Qui, {\it  Nucl.\ Phys.}\ {\bf B285}, 423 (1987);

\item{} A.\ Kato, {\it Mod.\ Phys.\ Lett.}\ {\bf A2}, 585 (1987);

\item{} P.\ Roberts, Ph.D.\ Dissertation, University of G\"oteborg, 1992;

\item{} P.\ Slodowy, {\it Bayreuther Math.\ Schr.}\ {\bf 33}, 197 (1990).

\item{3.} M.\ Hazewinkel, W.\ Hesselink, D.\ Siersma, and F.D.\ Veldkamp, 
{\it Nieuw Arch.\ Wisk.}\ {\bf 25}, 257 (1977);

\item{} P.\ Slodowy, in: {\it Algebraic Geometry}, Lecture Notes in Math
 1008, Springer, Berlin, 1983.

\item{4.} A.\ Ocneanu, Lectures at Fields Institute, April 26-30, 1995.

\item{5.} J.\ B\"ockenhauer and D.E.\ Evans, ``Modular invariants,
graphs and $\alpha$-induction for nets of subfactors II, III''
(preprints, hep-th/9805023, hep-th/9812110).

\item{6.} M.\ Bauer and C.\ Itzykson, {\it Commun.\ Math.\ Phys.}\ {\bf 127}, 
617 (1990).

\item{7.} T.\ Gannon, ``The level 2 and 3 modular invariants for the
orthogonal algebras'' (preprint, math.QA/9809020).
 
\item{8.} T.\ Gannon, ``Kac-Peterson, Perron-Frobenius, and
the classification of conformal field theories'' (preprint, q-alg/9510026);

\item{}  T.\ Gannon, ``The ${\cal A}_{\star}{\cal D}_{\star}{\cal
E}_7$-type invariants of affine algebras'' (in preparation).

\item{9.} W.\ Nahm, {\it Duke Math.\ J.}\ {\bf 54}, 579 (1987);

\item{} W.\ Nahm, {\it Commun.\ Math.\ Phys.}\ {\bf 118}, 171 (1988).

\item{10.} J.-B.\ Zuber, in: {\it Proceedings of the XIth International
 Conference of Mathematical Physics},
 International Press, Boston, 1995.

\item{11.} R.\ Behrend. P.\ Pierce, V.\ Petkova and J.-B.\ Zuber,
``On the classification of bulk and boundary conformal field
theories''
(preprint, hep-th/9809097).

\item{12.} A.\ Hanany and Y.-H.\ He, ``Non-abelian finite gauge theories''
(preprint, hep-th/9811183).

\item{13.} A.\ Cappelli, {\it Phys.\ Lett.}\ {\bf B185}, 82 (1987).

\item{14.} {P.\ Ginsparg}, {\it Nucl.\ Phys.}\ {\bf B295 [FS21]},
153 (1988);

\item{} {E.\ Kiritsis}, {\it Phys.\ Lett.}\ {\bf B217}, 427 (1989).

\item{15.} T.\ Gannon,  
{\it Nucl.\ Phys.}\ {\bf B491}, 659  (1997). 

\end